\documentclass[12pt]{article}
\usepackage{amssymb}

\def\header{A. Miller \hfill }
\def\col{{\mathbb Col}}

\def\forces{{\;\Vdash}}
\def\group{{\mathcal H}}
\def\la{\langle}
\def\name#1{\stackrel{\circ}{#1}}
\def\nn{{\mathcal N}}
\def\om{\omega}
\def\proof{\par\noindent Proof\par\noindent}
\def\qed{\par\noindent QED\par\bigskip}
\def\ra{\rangle}
\def\res{\upharpoonright}
\def\rmand{\mbox{ and }}
\def\rmor{\mbox{ or }}
\def\sm{{\setminus}}
\def\st{\;:\;} 
\def\su{\subseteq}

\newtheorem{theorem}{Theorem}
\newtheorem{prop}[theorem]{Proposition}

\newtheorem{define}[theorem]{Definition}

\begin{document}

\begin{center}
{\large Remark 3.4\\ A Dedekind Finite Borel Set}
\end{center}

\begin{flushright}
Arnold W. Miller\\
August 2012
\end{flushright}

\def\address{\begin{flushleft}
Arnold W. Miller \\
miller@math.wisc.edu \\
http://www.math.wisc.edu/$\sim$miller\\
University of Wisconsin-Madison \\
Department of Mathematics, Van Vleck Hall \\
480 Lincoln Drive \\
Madison, Wisconsin 53706-1388 \\
\end{flushleft}}

Asaf Karagila pointed out that Remark 3.4 in my paper,
A Dedekind finite Borel set \cite{ded},
directly contradicts Theorem 3.3 (c) in
Does GCH imply AC locally? by A. Kanamori and D. Pincus \cite{kanamori}.
Neither of the papers provides a
proof.\footnote{
Aki Kanomori tells me that
Theorem 3.3 (c) was incorrectly stated.
Instead of $\nn$ it
should be the Solovay-type inner model $\nn_0\subseteq \nn$.
This is the model
$HOD(^\omega ON)$ of $\nn$.
See page 155 of {\bf Consequences of the
Axiom of Choice} by Rubin and Howard who call
$\nn_0$ Truss's model.}
Here we give a proof of Remark 3.4 \cite{ded}.

The Feferman-Levy model $\nn$ is described in
Cohen \cite{cohen} p.143, Jech
\cite{jech} p.142,
and in the on-line only appendix to Miller \cite{longbor}.
We use the notation from \cite{longbor}.
$\col$ is the Levy collapse of $\aleph_\om$ , $\col_n$
the collapse of $\aleph_k$ for $k<n$, $G_n=\col_n\cap G$
where $G$ $\col$-generic
over $M$. The model $\nn$ is the symmetry model
$M\su\nn\su M[G]$ determined by the group $\group$ of
automorphism of $\col$ which finitely permute finitely
many of the collapsing
maps domains and the filter of subgroups generated by
$(H_n\;:\;n<\om)$ where
$$H_n=\{\rho\in \group\st \forall p\in\col_n\;\;
\rho(p)=p\}.$$

\begin{define}
Let $\la,\ra:\om\times\om\to\om$ be a fixed bijection,
i.e., a pairing function.  For each $n\in\om$
define the map $\pi_n:2^\om\to 2^\om$ by:
$$\pi_n(x)=y \mbox { iff } \forall m\in\om\;\;y(m)=x(\la n,m\ra).$$
\end{define}

\begin{define}
In the Feferman-Levy model $\nn$, take $F_n=M[G_n]\cap 2^\om$.
Define
$$B=\{x\in 2^\om\st\forall n\;\; \pi_n(x)\in F_{n+1}\sm F_n
\rmor [\pi_n(x)\in F_n \rmand \pi_n(x)=\pi_{n+1}(x)]
\}.$$
\end{define}

The set $B$ is uncountable because there is a map $h$ from $B$ onto $2^\om$.
Define $h$ by $h(x)=\pi_n(x)$ iff $\pi_n(x)=\pi_m(x)$ for all $m>n$.  Such an
$n$ must exist because for any $x$ there exists $n$ such that $x\in F_n$ and
hence $\pi_m(x)\in F_n$ for all $m$.  It is easy to check that $h$ maps $B$ onto
$2^\om$.

\begin{prop}
In the Feferman-Levy model $\nn$ the set $B$ has the
property that there is no one-to-one map
taking $2^\om$ into $B$.
\end{prop}

\proof
Suppose for contradiction that
there is in $\nn$ a one-to-one map  $f:2^\om\to B$.
Define
$$C_n=\{x\in B\st \forall m>n \;\;\pi_m(x)=\pi_n(x)\}$$
and note that $B=\bigcup_nC_n$ and each $C_n$ is countable.

Let
$$p\forces \name{f}:2^\om\to \name{B} \mbox{ is one-to-one }$$
and suppose that $p\in\col_n$ and $H_n$
(the subgroup of $\group$ which
is the identity on $\col_n$)
fixes $\name{f}$.  Since
$C_n$ is countable and $f$ is one-to-one, we may choose
$N>>n$ and $q\leq p$ such that
$$q\forces {f}^{-1}(\name{C}_n)\su \name{F}_N.$$
Now let $\name{x_N}$ be a canonical name for a real which codes
the generic collapse of $\aleph_N$.  The important property
it has is that $x_N\in F_{N+1}\sm F_{N}$ and its values are
completely determined by $N^{th}$ coordinate of the generic
filter and not any earlier
ones, e.g., the $\leq n$ coordinates. Take $r\leq q$ such that
$$r\forces \name{f}(\name{x}_N)=\name{y}, \pi_n(\name{y})=\tau, \rmand
\tau\in \name{F}_{n+1}\sm \name{F}_n$$
where $\tau$ is a $\col_{n+1}$-name for an element of $2^\om$.
For $G$ $\col$-generic over $M$ and containing
$r$ note that $\tau^{G_{n+1}}\in M[G_{n+1}]\sm M[G_{n}]$.  So we can
find $k$ and $r_0,r_1\leq r$ such that
\begin{itemize}
\item $r_0$ and $r_1$ are the same on every coordinate except $n$,
\item $r_0\forces \tau(k)=0$, and
\item $r_1\forces \tau(k)=1$.
\end{itemize}

Now let $\rho\in\group$ be an automorphism of $\col$
which fixes every coordinate except
$n$ and makes $\rho(r_0)$ and $r_1$ compatible.
To do this move the domain of $r_0\res(\{n\}\times\om)$ disjoint from the
domain of $r_1$.
Since $\rho\in H_n$ it fixes $\name{f}$, i.e.,
$\rho(\name{f})=\name{f}$.   Since $\rho$ fixes the $N^{th}$
coordinate $\rho(\name{x}_N)=\name{x}_N$.
But this means that
$$\rho(r_0)\forces \pi_n(\name{f}(\name{x}_N))(k)=0$$
and
$$r_1\forces \pi_n(\name{f}(\name{x}_N))(k)=1$$
which contradicts that $\rho(r_0)$ and $r_1$ are compatible.
\qed

\address


\begin{thebibliography}{99}

\bibitem{cohen} Cohen, Paul J.; {\bf Set theory and the
continuum hypothesis.} W. A. Benjamin, Inc.,
New York-Amsterdam 1966 vi+154 pp.

\bibitem{jech} Jech, Thomas J.; {\bf The axiom of choice.} Studies in Logic and
the Foundations of Mathematics, Vol. 75. North-Holland Publishing Co.,
Amsterdam-London; American Elsevier Publishing Co., Inc., New York, 1973. xi+202
pp.

\bibitem{kanamori}
Kanamori, A.; Pincus, D.;
Does GCH imply AC locally?,
Paul Erdos and his mathematics, II (Budapest, 1999), 413-426,
Bolyai Soc. Math. Stud., 11, János Bolyai Math. Soc., Budapest, 2002.

http://math.bu.edu/people/aki/7.pdf


\bibitem{longbor} Miller, Arnold W.;  Long Borel hierarchies, Math Logic
Quarterly, 54(2008), 301-316.

http://www.math.wisc.edu/$\sim$miller/res/longbor.pdf

\bibitem{ded} Miller, Arnold W.;
A Dedekind Finite Borel Set,
Arch. Math. Logic 50 (2011), no. 1-2, 1--17.

http://www.math.wisc.edu/$\sim$miller/res/ded.pdf

\end{thebibliography}
\end{document}